\documentclass[a4paper, 12pt]{article}
\usepackage[utf8]{inputenc}
\usepackage[T2A]{fontenc}
\usepackage{amsmath,amssymb,amsthm, amsfonts, mathtools, bm}
\usepackage[a4paper,hmargin=2.5cm,vmargin=2.5cm]{geometry}
\usepackage[english]{babel}

\newtheorem{thm}{Theorem}

\begin{document}

\makeatletter

\title{Universal matrix Capelli identity}
\author{\rule{0pt}{7mm} Mikhail Zaitsev \thanks{mrzaytsev@edu.hse.ru}\\
{\small\it
National Research University Higher School of Economics}\\
}

\maketitle

\begin{abstract}
We propose a universal matrix Capelli identity and explain how to derive Capelli identities for all quantum immanants in the Reflection Equation algebra and in the universal enveloping 
algebra $U(\mathbf{gl}_{M|N})$.
\end{abstract}


{\bf Keywords:} Capelli's identity, Reflection Equation algebra, quantum immanants

\vspace{\baselineskip}

\noindent
{\bf 1.} We present a universal matrix Capelli identity in the Reflection Equation algebra (Theorem 1). It is a far-reaching generalization of the classical Capelli identity in the Weyl algebra (see [C]):
\begin{equation} \label{class_cap}
    cdet(XD + K) = detX detD.
\end{equation}
Here $X = || x_i^j ||_{1\leq i,j \leq N}$ is a matrix of commuting variables, $D = || \frac{\partial}{\partial x_j^i}||_{1\leq i,j \leq N}$ is a matrix of partial derivatives and $K = diag(N-1,N-2, \dots,1,0)$ is a diagonal numerical matrix. The symbol $cdet$ denotes the column determinant of a non-commutative matrix.

There are known many generalizations of this identity. In particular, the paper \cite{OK1} proposes a version of the identity for quantum $\mu$-immanants associated with Young diagrams. For example, the identity (\ref{class_cap}) corresponds to a single-column diagram with $N$ cells. In \cite{OK2} a matrix version of this identity is proposed; the previous version is obtained by taking traces at all positions. This version is called matrix, as well as all versions in which traces are not taken.

Further generalizations of the Capelli identity are related to the Reflection Equation algebras, corresponding to  the Hecke type $R$-matrices,  they are called {\em quantum identities}. This name is motivated by the fact that if $R$ is the Drinfeld-Jimbo $R$-matrix, the corresponding Reflection Equation algebra is a deformation (quantization) of the universal enveloping algebra $U(\mathbf{gl}_N)$. Quantum matrix Capelli identities for the single-row and single-column diagrams were obtained in \cite{GPS1} (for any skew-invertible $R$-matrix).
 In \cite{JLM} this result was  generalized to other Young diagrams but only for  the  Reflection Equation algebras
 associated with the Drinfeld-Jimbo R-matrices. In this case quantum analogs of  the immanants, introduced by A. Okounkov, also have been introduced in \cite{JLM}. Note that the universal matrix Capelli identity (\ref{Cap-quant}) holds in  the Reflection Equation algebras corresponding  to all skew-invertible Hecke $R$-matrices and it is not attached to
 any Young diagram. All other forms of Capelli identity can be  easily deduced from it.

\medskip

\noindent
{\bf 2.} For any matrix $A \in Mat_{N \times N}( \mathbb{C}) \otimes U$ with coefficients in an arbitrary vector space $U$ we use the notation:
$$
A_i = \mathbb{ E} \otimes \mathbb{ E} \otimes \ldots  \underset{i}{A} \ldots  \otimes \mathbb{ E} \in Mat_{N \times N}( \mathbb{C})^{\otimes k}\otimes U,
$$
where $\mathbb{E}$ is the identity matrix in the corresponding tensor component. Using this notation, the relations in the universal enveloping algebra $U(\mathbf{gl}_N)$ can be written in the following form:
$$L_1 L_2 - L_2 L_1 = L_1 P_{12} - P_{12} L_1,$$
where $L= ||L_i^j||_{1 \leq i,j \leq N}$ denotes the matrix composed of the algebra generators and $ P_{12}$ is a permutation matrix with the matrix elements $P_{kl}^{mn}=\delta_k^n\, \delta_l^m$, $ 1\leq k,l,m,n\leq N$.

Let us consider the Weyl algebra $W_N$ generated by the elements of the matrices $X$ and $D$. These matrices satisfy the following relations:
$$
X_1 X_2 = X_2 X_1, \quad D_1 D_2 = D_2 D_1,\quad  D_1 X_2 = X_2 D_1 + P_{12}.
$$
The algebra $U(\mathbf{gl}_N)$ is embedded into the Weyl algebra via an embedding that maps the elements of the matrix $L$ to the corresponding elements of the matrix $XD$. The universal matrix Capelli identity related to the algebra $U(\mathbf{gl}_N)$  can be written in the following form:
\begin{equation}
L_{1} (L_{2} - j_2) \dots (L_{n} - j_n) = X_{1} \dots X_{n} D_{1} \dots D_{n},\quad \forall\,n\ge 1,
\label{Cap-cl}
\end{equation}
where the elements $j_k = \sum\limits_{i = 1}^{k-1} P_{i k}$ are the images of the Jucys-Murphy elements defined  in the group algebra  $\mathbb{C}[S_n]$    of the symmetric group.
As a straightforward  consequence one can get  the Capelli identity from \cite{OK1}:
\begin{equation}
Tr_{(1 \dots n)} \left(L_1 (L_2 - c(2)) \dots (L_n - c(n)) P_{ii}^{\lambda} \right) = Tr_{(1 \dots n)} \left(X_1 \dots X_n D_1 \dots D_n P_{ii}^{\lambda} \right).
\label{imm}
\end{equation}
Here we use the notation $P_{ii}^{\lambda}$ for the primitive idempotents of the group algebra $\mathbb{C}[S_n]$. Observe that the primitive idempotents are labeled by the standard Young tableaux
$(\lambda, i)$, where the partition $\lambda\vdash n$ defines the form of the tableau and index $i$ numerates the  tableaux  of shape $\lambda$. The symbol $c(k)$ denotes the content of cell with number $k$ in the given standard tableau. We also use the notation $Tr_{(k)}$ for the trace in the $k$-th tensor component.
Note that both parts of the equality (\ref{imm}) do not depend on the number $i$ of the primitive idempotent $P_{ii}^{\lambda}$. Besides, on omitting the trace in (\ref{imm}), we get the matrix Capelli identity from \cite{OK2}.

\medskip

\noindent
{\bf 3.} The Reflection Equation algebra $M(R)$ is a unital  associative algebra generated by $N^2$ generators $m_i^j$ subject to
\begin{equation}
R_{12}M_1R_{12}M_1 = M_1R_{12}M_1R_{12},\quad  M = \| m_i^j \|_{1\leq i,j \leq N}.
\label{REA-M}
\end{equation}
Here $R \in Mat_{N \times N}(\mathbb{C})^{\otimes 2}$ is the skew-invertible Hecke $R$-matrix, i.e. the following properties are satisfied:
$$
R_1R_2R_1 = R_2R_1R_2, \qquad R^2 = 1 + (q-q^{-1}) R, \quad q \in \mathbb{C}\backslash \{ 0, \pm 1\}
$$
and there exists a matrix $\Psi$ such that $Tr_{(2)}(R_{12} \Psi_{23}) = P_{13} = Tr_{(2)}( \Psi_{12}R_{23})$.
Here and below we use notation $R_i = R_{ii+1}$ and assume that the parameter $q$ is generic, that is, $q^{k} \neq 1,$ $\forall k \in \mathbb{Z}_{>0}$. Using the linear change of the 
generators $\hat L = \frac{{\Bbb E} - M}{q - q^{-1}}$ we reduce the relation (\ref{REA-M}) to the form:
\begin{equation}
		\hat L_1R_{1}\hat L_1R_{1} - R_{1}\hat L_{1}R_{1}\hat L_{1} = \hat L_{1} R_{1} - R_{1} \hat L_{1}.
		\label{mREA}
		\end{equation}
In particular, for the Drinfeld-Jimbo $R$-matrix, the algebra defined by (\ref{mREA}) is a deformation of the universal enveloping algebra $U(\mathbf{gl}_N)$. Moreover, in the case when a Hecke $R$-matrix is a deformation of the super-permutation $P_{M|N}$, this algebra is a deformation of the algebra $U(\mathbf{gl}_{(M|N)})$.

Let us introduce a convenient notation: for an arbitrary $N\times N$ matrix $A$ we set
$$
A_{\overline{k}} = R_{k-1} \dots R_{1} A_1 R_{1}^{-1} \dots R_{k-1}^{-1}.
$$
Analogically to the classical case, the algebra $M(R)$ is embedded in the quantum Weyl algebra $W(R)$, which is  the quotient-algebra of the free product of two Reflection Equation algebras $M(R)$ and $D(R^{-1})$ (in the definition of the latter algebra, $R$ is replaced by $R^{-1}$). The generating matrices of these two algebras are subject
to the following permutation relations:
$$
D_1 M_{\overline{2}} = R_1^{-1} + M_{\overline{2}}D_1 R_{1}^{-2}.
$$
(These relations define the corresponding ideal.)
It is known (see \cite{GPS1}) that the matrix $\hat L=MD$ satisfies the relation (\ref{mREA}). Thus,  the Reflection Equation algebra is represented by quantum analogs of vector fields. Now we can formulate the main theorem.
\begin{thm}
     For any integer $n\ge 1$  the universal quantum matrix Capelli identity holds:
				\begin{equation} \label{Cap-quant}
                  \hat L_{\overline{1}} \left(\hat L_{\overline{2}} + \frac{J_{2}^{-1} -1}{q-q^{-1}}\right) \dots \left(\hat L_{\overline{n}} + \frac{J_{n}^{-1} -1}{q-q^{-1}}\right) =
									M_{\overline{1}} \dots M_{\overline{n}} D_{\overline{n}} \dots D_{\overline{1}} J_{1}^{-1} \dots J_{n}^{-1}.
				\end{equation}
        Here $J_1 = 1$, $J_k = R_{k-1} \dots R_2 R_{1}^2 R_2 \dots R_{k-1}$, $k >1$ are the images of the Jucys-Murphy elements of the Hecke algebra with respect to the $R$-matrix representation.
    \end{thm}
Note that the word ``quantum'' means that the identity is related to the Reflection Equation algebra. If $R$-matrix is a deformation of the permutation matrix, the universal matrix Capelli identity (\ref{Cap-cl}) in the algebra $U(\mathbf{gl}_{N})$ is obtained from (\ref{Cap-quant}) by passing to the limit $q \rightarrow 1$.

Analogously, by taking any Hecke $R$-matrix, that deforms the super-permutation $P_{M|N}$, and by passing to the limit $q \rightarrow 1$, we obtain a similar identity in the superalgebra $U(\mathbf{gl}_{(M|N)})$.

\medskip

\noindent
{\bf 4.} Note that for generic $q$ the Hecke algebra $\mathbb{H}_n(q)$ of type $A_{n-1}$ is semisimple and is isomorphic to the group algebra of the symmetric group ${\Bbb C}[S_n]$.
In particular, the primitive idempotents $E_{ii}^{\lambda}$ of the algebra $\mathbb{H}_n(q)$ are labeled by the standard Young tableaux. The set of Jucys-Murphy elements generates a maximal commutative subalgebra of the Hecke algebra and their action on the primitive idempotents is given by  $J_k E_{ii}^{\lambda} = q^{2c(k)} E_{ii}^{\lambda}$.
Multiplying the equality (\ref{Cap-quant}) by the idempotent $E_{ii}^{\lambda}$ and using the last relation, we obtain the following identity:
\begin{equation}
\hat L_{\overline{1}} (\hat L_{\overline{2}} - q^{-c(2)}[c(2)]_{q}) \dots (\hat L_{\overline{n}} - q^{-c(n)}[c(n)]_{q}) E_{ii}^{\lambda} =
q^{-2(c(1) + \dots + c(n))}M_{\overline{1}} \dots M_{\overline{n}} D_{\overline{n}} \dots D_{\overline{1}} E_{ii}^{\lambda}.
\label{corCap}
\end{equation}
The quantum matrix Capelli identity from \cite{JLM} was proved in this form. Also, in this paper   quantum immanants were introduced in the Reflection Equation algebra associated with the Drinfeld-Jimbo $R$-matrix.
In the general case quantum immanants can be obtained from the left-hand side of the equality (\ref{corCap}) by applying the $R$-traces $Tr_{R(12...n)}$ in all tensor components. This definition does not appeal to the form of the $R$-matrix.

\end{document}